\documentclass[11pt]{amsart}
\usepackage{epsfig}

\newcommand{\op}{\mbox{op}}

\newtheorem{theorem}{Theorem}[section]
\newtheorem{lemma}[theorem]{Lemma}
\newtheorem{proposition}[theorem]{Proposition}

\theoremstyle{definition}

\newtheorem{example}{Example}

\newtheorem{cor}{Corollary}

\theoremstyle{remark}

\numberwithin{equation}{section}

\title[Quantum knots invariants]{On the holomorphic point of view in the theory
of quantum knot invariants}
\author{R{\u{a}}zvan Gelca}
\address{Department of Mathematics and Statistics, 
Texas Tech University, Lubbock, TX 79409 and Institute of Mathematics
of the Romanian Academy, Bucharest, Romania}
\email{rgelca@gmail.com}
\thanks{This paper appeared in the Journal of Geometry and Physics,
56 (2006).}
\subjclass{81T45, 57M27, 81R50, 81S10}

\date{9 March 2005}

\keywords{Witten-Reshetikhin-Turaev invariants,  
theta functions, Weyl quantization, modular functor}

\begin{document}
\maketitle

\begin{abstract}
In this paper we describe progress made toward
the construction of the Witten-Reshetikhin-Turaev theory
of knot invariants from a geometric point of view. This
is done in the perspective of a joint result of the author
with A. Uribe which relates the quantum group 
and  the Weyl quantizations 
of the moduli space of flat $SU(2)$-connections on the
torus. Two results are emphasized: the reconstruction
from Weyl quantization of the restriction to the torus 
of the modular functor, 
and a description of a basis of the space
of quantum observables on the torus 
in terms of colored curves, which answers a question
related to quantum computing.  

\end{abstract}

\section{Introduction}

It is known that, for a compact  simple Lie group $G$, 
  the Hilbert space of the quantization of the 
 moduli space of flat $G$-connections on a surface is the space
of holomorphic sections of the Chern-Simons line bundle. Alternately,
this space has a basis consisting of 
 admissible colorings  of  
 the core of the handlebody bounded by the surface 
by irreducible representations of the quantum
group of $G$.

The quantization of the moduli space of flat $SU(2)$-connections on
the {\em torus} was studied in detail in \cite{GU}. In that work
two quantization models were compared, the equivariant Weyl quantization 
of the complex plane that covers the
moduli space and the quantum group quantization,  performed 
after Witten's ideas \cite{Wi} with the techniques of
Reshetikhin and Turaev \cite{RT}, \cite{T2}. 
There, it was shown that these two 
quantizations are {\em unitarily equivalent}.

This result gives rise to new possibilities for developing 
$SU(2)$ Chern-Simons theory, and in particular the study of 
the Jones polynomial \cite{Jones}, from a
geometric point of view.  
The present paper describes some progress made in this direction.
For other contributions to the subject see
 \cite{jea}, \cite{jeau1}, \cite{jeau2}, \cite{freed}.

Here is a description of the contents of the paper.
Section 2 gives a brief overview of the two quantization models. 
In section 3 we show how the projective representation of the mapping
class group of the torus, which arises in the Reshetikhin-Turaev 
topological quantum field 
theory, can be recovered from the Weyl quantization model. It is 
important to note that the Weyl  quantization  contains all the 
necessary information
about the modular functor restricted to the torus. Next 
section contains  explicit descriptions of quantum knot 
invariants as holomorphic sections of the Chern-Simons line bundle.
Section 5 is devoted to the quantum observables. They are 
described as integral operators, much in the spirit of Witten's path
integral, then their spectra are computed. The paper ends
with an application to quantum computing. It determines the 
 basis of the Hilbert space  of ground states of a certain quantum system,
which is this same as the vector space of the 
Chern-Simons theory of  the quantum double of the group $SU(2)$. 
This solves  a problem that arose in \cite{fnsww}.
The author would like to thank the referee, whose comments
considerably improved the quality of the manuscript.

\section{The quantization of the moduli space of flat
$SU(2)$-connections on the torus}

The Hilbert space of the
quantization of the moduli space of flat $SU(2)$-connections 
on the torus was described in \cite{GU}. Since in that paper
 a small error occurred
in the  exposition, but fortunately not in the final result, we briefly 
explain the construction again. 

Denote by ${\mathcal M}$  
the moduli space of flat $SU(2)$-connections on the
torus ${\mathbb T}^2$,
 which is known informally as the ``pillow case''. ${\mathcal M}$ 
is the quotient of the complex plane by the symmetries $z\rightarrow
z+m+ni$, $m,n\in {\mathbb Z}$, and $\sigma(z)= -z$. The symplectic
form  which determines the Poisson bracket and
consequently the ``classical mechanics'' on ${\mathcal M}$ is
 $\omega = -{\pi}dz\wedge d\bar{z}=2\pi idx\wedge dy$. 

The Hilbert space of the quantization consists of holomorphic 
 $1$-forms on the smooth part of ${\mathcal M}$ 
with values in a line bundle of curvature 
$\omega ^{N}$, where $N=1/\hbar$ is the reciprocal of Planck's constant.
Constrains given by the Weil integrality condition and the
Reshetikhin-Turaev theory require $N$ to
be an even integer, $N=2r$. In Witten's theory, the number $r-2$ is
called the level of the quantization.  
The inner product is defined by integrating  the cup product of
two 1-forms over ${\mathcal M}$. 

A $1$-form on ${\mathcal M}$ can be written locally as
$f(z)dz$, where $f(z)$ is a section of the line bundle. This could almost
be done globally, by lifting it to ${\mathbb C}$, except that on ${\mathbb C}$
the symmetry with respect to the origin $\sigma$
 changes $dz$ to $-dz$. This sign
change can be incorporated into the line bundle. Consequently, the line bundle 
is  determined by a cocycle $\chi$ on ${\mathbb C}$ satisfying
\begin{alignat*}{1}
\chi (z, m+in) &= \exp 2r\pi(-2inz+n^2), 
\quad m, n \in {\mathbb Z},\\
 \chi (z, \sigma) &=-1.
\end{alignat*} 

The  Hilbert space of the quantization is then  identified with the space of
odd theta functions 
\begin{eqnarray*}
{\mathcal H}_r=\{f\in{\mathcal Hol}({\mathbb C})|
 \, f(z+m+in)=e^{{2r}\pi(n^2-2inz)}f(z),
 \, f(z)=-f(-z)\}.
\end{eqnarray*}

If we denote
\begin{eqnarray*}
\theta_j(z)=
\sum_{n=-\infty}^{\infty}
e^{-\pi(2rn^2+2jn)+2\pi iz(j+2rn)},
\end{eqnarray*}
then the functions 
\begin{eqnarray*}
\zeta_j(z)=
\sqrt[4]{r}e^{-\pi j^2/2r}(\theta_j(z)-\theta_{-j}(z)), \quad j=1,2,\ldots, r-1
\end{eqnarray*}
form an orthonormal basis of ${\mathcal H}_r$. For later use we extend
the definition of $\zeta_j(z)$ to all integers $j$ by the relations
$\zeta_{r+j}(z)=-\zeta_{r-j}(z)$ and
$\zeta_{-j}(z)=-\zeta_{j}(z)$. Note that  $\zeta_0(z)=\zeta_r(z)=0$. 

The quantum observables are defined using Weyl quantization. 
For a classical observable $f\in C^{\infty}({\mathcal M})$, 
 the associated quantum observable is the Toeplitz
operator with symbol $e^{-\frac{\Delta}{8r}}f$. Here $\Delta $
is the Laplace operator
$\frac{1}{2\pi}\left(\frac{\partial^2}{\partial x^2}
+\frac{\partial ^2}{\partial y^2}\right)$. We use view  Weyl
quantization 
as interpolating between Wick and anti-Wick
as explained in \cite{folland}.

The alternative quantization model, which
appears in the realm of quantum groups \cite{RT}, is the following.
The  orthonormal basis of the Hilbert space is formally
identified with  $V^j( \alpha)$, $j=1,2,\ldots, r-1$,  the colorings
 of the 
core $\alpha $ of the solid torus  by the $j$-dimensional
irreducible representations of the quantum group of 
$SL(2,{\mathbb C})$ at a root of unity $q=\exp(2\pi i/r)$. 
The algebra of classical observables contains as a dense subset the
ring generated by the traces of holonomies of $SU(2)$-connections along
simple closed curves on the torus. 
The quantum observable associated to the trace in the $k$-dimensional 
irreducible representation of $SU(2)$ of such a  holonomy along a curve 
$\gamma$ is
simply the curve $\gamma$ colored by the $k$-dimensional 
irreducible representation
of the quantum group of $SL(2,{\mathbb C})$ (with the usual convention
when $k$ is larger than $r$, the level of the quantization). 
There is a way of identifying colored curves with  linear operators using
knot invariants. For details we refer the reader to \cite{GU}.

The main result of \cite{GU} is that the unitary
map $\zeta_j(z)\rightarrow V^j(\alpha)$, $j=1,2,\ldots, r-1$, is an
equivalence between the Weyl quantization 
and the quantum group quantization. In the process of proving this result
we discovered a formula which will be used in the sequel. 
  Denote by $C(p,q)$, $p,q\in {\mathbb Z}$ 
the operator obtained by performing Weyl quantization with
 symbol $2\cos 2\pi (px+qy)$. Then
\begin{eqnarray*}
C(p,q)\zeta_m(z)=t^{-pq}(t^{2qm}\zeta_{m-p}(z)+t^{-2qm}\zeta_{m+p}(z)), \quad 
m=1,2,\ldots, r-1,
\end{eqnarray*}
where $t=\exp(\pi i/2r)$ (and so $t^4=q$). 

This paper describes some features of the 
Witten-Reshetikhin-Turaev theory from the analytical-geometric point of view. 
We will  need  an alternative formula for $\zeta _j(z)$.
 To obtain it we write $\zeta_j(z)$ as
\begin{eqnarray*}
\sqrt[4]{r}
e^{-\frac{\pi j^2}{2r}}\sum _{n=-\infty}^{\infty}\left(
e^{-\pi (2rn^2+2jn)+2\pi i z(j+2rn)}-e^{-\pi (2rn^2-2jn)+2\pi iz(-j+2rn)}
\right).
\end{eqnarray*}
Note that 
\begin{alignat*}{1}
e^{-\frac{\pi j^2}{2r}}
e^{-\pi (2rn^2+2jn)+2\pi i z(j+2rn)}&=
 e^{-\frac{\pi j^2}{2r}}
e^{2\pi izj}e^{-\pi (2rn^2+2jn)+4\pi irnz}\\
&= e^{-2\pi rz^2} e^{-2\pi r
\left(n+\frac{j}{2r}-iz\right)^2}.
\end{alignat*}
Using the Poisson formula for $e^{-x^2}$, the first sum can be  transformed 
into
\begin{eqnarray*}
 e^{-2\pi rz^2}\sum _{n=-\infty}^\infty
e^{-\frac{\pi n^2}{2r}}e^{2\pi i \left(\frac{j}{2r}-iz\right)n}
= \sum _{n=-\infty}^\infty e^{-2\pi rz^2-\frac{\pi n^2}{2r}+2\pi nz}
e^{\frac{\pi ijn}{r}} .
\end{eqnarray*}
Replacing $j$ by $-j$ we obtain  that the other sum is equal to 
\begin{eqnarray*}
\sum _{n=-\infty}^\infty e^{-2\pi rz^2-\frac{\pi n^2}{2r}+2\pi nz}
e^{-\frac{\pi ijn}{r}}.
\end{eqnarray*}
Subtracting the two and recalling the definition of  the quantized integer 
\begin{eqnarray*}
[n]=\frac{e^{\frac{\pi in}{r}}-e^{-\frac{\pi in}{r}}}{e^{\frac{\pi i}{r}}-
e^{-\frac{\pi i}{r}}}=\frac{\sin\frac{n\pi}{r}}{\sin \frac{\pi}{r}}=
\frac{t^{2n}-t^{-2n}}{t^2-t^{-2}}. 
\end{eqnarray*}
we obtain  the following:

 \begin{lemma}\label{formula}
For $j=1,2,\ldots, r-1$, 
\begin{eqnarray*}
\zeta_j(z)=2i\sqrt[4]{r}\sin\frac{\pi}{r}
\sum_{n=-\infty}^\infty
e^{-2\pi r\left( z-\frac{n}{2r}\right)^2}[nj].
\end{eqnarray*}
\end{lemma}

As a corollary  we find that  that the reproducing kernel
of the  Hilbert space ${\mathcal H}_r$ is 
\begin{alignat*}{1}
K(z,w)&=2r^{\frac{3}{2}}\sum_{\begin{array}{l}
n,m=-\infty\\ 2r|m-n 
\end{array}}^\infty e^{-2r\pi\left[\left(z-\frac{n}{2r}\right)^2+
\left(\overline{w}-\frac{m}{2r}\right)^2\right]}\\
&-
2r^{\frac{3}{2}}\sum_{\begin{array}{l}
n,m=-\infty\\ 2r|m+n
\end{array}}^\infty e^{-2r\pi\left[\left(z-\frac{n}{2r}\right)^2+
\left(\overline{w}-\frac{m}{2r}\right)^2\right]}\\
&=4r^{\frac{3}{2}}\sum_{k,m=-\infty}^\infty e^{-2r\pi \left[
(z+k)^2+\frac{m^2}{4r^2}+\left(\overline{w}-\frac{m}{2r}\right)^2\right]}
\sinh m\pi (z+k).
\end{alignat*}
This means that if $f(z)$ is an element of the Hilbert space of the
quantization, namely an odd theta function,
 then $f(z)=\int_{\mathcal M}f(w)K(z,w)dw$, for all $z\in
{\mathcal M}$. 

\section{The projective representation of the mapping class group of the
torus on the Hilbert space of the quantization}

Part of the Reshetikhin-Turaev topological quantum field theory
is a projective representation of the mapping class group of the torus
onto the Hilbert space of the quantization. This representation
is the restriction of
the modular functor to the torus.
We will show how this
projective representation can be recovered from Weyl quantization.

A simple closed curve $\gamma$ on the torus defines a smooth function
on the moduli space of flat $SU(2)$-connections by taking the 
trace in the fundamental representation of $SU(2)$ of the holonomy along
$\gamma $ of the connection. Call this function $f_{\gamma}$ and
$\mbox{op}(f_{\gamma})$ the operator associated to it through Weyl quantization.

There is a right action of the mapping class group of the torus on
$C^{\infty}({\mathcal M})$. It is defined as follows.
Identify the moduli space  ${\mathcal M}$  with the
algebraic variety of characters of  $SU(2)$-representations of the fundamental
 group of the torus. 
If $g$ is an element of the mapping class group of the torus 
and $\sigma$ is an $SU(2)$-representation of the fundamental group of
the torus, then
\begin{eqnarray*}
g\cdot \sigma :\quad \gamma \rightarrow \sigma(g(\gamma)),
\end{eqnarray*}
We therefore have a 
left action of the mapping class group of the torus on the moduli
space ${\mathcal M}$, which then induces a right action on $C^{\infty}(
{\mathcal M})$ by $(f\cdot g)(x)=f(g\cdot x)$. Through Weyl quantization we
obtain a right action 
on quantum observables.

In particular, if the element $g$ of the mapping class
group maps the curve $\gamma$ to $g(\gamma)$, then $f_{\gamma}\cdot g=
f_{g(\gamma)}$. In this sense we have a natural action of the mapping
class group on symbols of operators and therefore a natural action
on operators themselves.

If we view the torus as ${\mathbb R}^2/{\mathbb Z}^2$, then
its mapping class group is generated by the maps ${\mathcal S}$ and $
{\mathcal T}$ defined
by ${\mathcal S}(x,y)=(-y,x)$, 
${\mathcal T}(x,y)=(x,x+y)$.  They act on  quantum observables by 
\begin{eqnarray*}
\quad  \op (f(x,y))\cdot {\mathcal S}=\op (f(-y,x)),\, 
\op(f(x,y))\cdot {\mathcal T}=\op (f(x,x+y)).
\end{eqnarray*}

The Reshetikhin-Turaev topological quantum field theory comes with a projective
representation of the mapping class group of the torus on ${\mathcal H}_r$
 defined by 
 $\rho({\mathcal S})= S$ and $\rho({\mathcal T})=T$, where
\begin{eqnarray*}
S=
\left([jk]\right)_{1\leq j,k\leq r-1}, 
\quad \mbox{and} \quad T=(\delta_{j,k}t^{j^2-1})_{1\leq j,k\leq r-1}.
\end{eqnarray*}
Here, as before,   $[jk]$ is the quantized integer, $t=e^{\frac{i\pi}{2r}}$, 
 and $\delta_{j,k}$ is the Kronecker symbol (as we are only
 interested
in projective representations we did not incorporate the factor
$1/X$, with  $X=\sqrt{\sum_{k=1}^{r-1}[k]^2}$ in the definition of
 $S$).

In order for the entire theory to be consistent, this representation
must be compatible with the natural action on the algebra of quantum 
observables, which means  that if $op(f)$ is a quantum observable and $g$ and
element of the mapping class group, then
\begin{eqnarray*}
op(f)\cdot g= \rho(g)^{-1}op(f)\rho(g).
\end{eqnarray*}
The next result shows that Weyl quantization together with this condition
determine the projective representation of the mapping class group 

\begin{theorem}\label{unique}
There is a unique projective representation of the mapping
class group of the torus on the Hilbert space of the quantization
which is compatible with the natural action of the mapping
class group on quantum observables, and this is the projective
representation from the Reshetikhin-Turaev theory.
\end{theorem}

\begin{proof} Let $S=(a_{k,j})_{1\leq k,j\leq r-1}$ and $T=
(b_{k,j} )_{1\leq k,j\leq r-1}$ be the $S$ and $T$-matrices of
such a projective representation. 

First, let us extend the definition of $a_{k,j}$ so that the indices
can be any integer numbers. The equalities
\begin{eqnarray*}
&& \sum_{k=1}^{r-1}a_{k,2r-j}\zeta _k(z)=S\zeta _{2r-j}(z)
=S(-\zeta _j(z))=-S\zeta_j(z)=-\sum_{k=1}^{r-1}
a_{k,j}\zeta _k(z)\\
&& \sum_{k=1}^{r-1}a_{k+r,j}\zeta_{k+r}(z)=S\zeta _j(z)=
\sum _{k=1}^{r-1}a_{k,j}\zeta _k(z)
\end{eqnarray*}
show that the correct choice for $a{k,j}$ is as odd and periodic in 
$j$ and $k$ of period $2r$, that is $a_{k,-j}=a_{-k,j}=-a_{k,j}$.
Moreover $\zeta _0=0$ means that we can choose $a_{k,0}=a_{0,j}=0$. 
Same conventions  for $b_{k,j}$. 

Recall that  $C(p,q)$, $p,q\in {\mathbb Z}$ denotes 
the operator with symbol $2\cos 2\pi (px+qy)$. We will use the 
previously mentioned formula 
\begin{eqnarray*}
C(p,q)\zeta_m(z)=t^{-pq}(t^{2qm}\zeta_{m-p}(z)+t^{-2qm}\zeta_{m+p}(z)).
\end{eqnarray*}

Let us look at the action of $S$ on
the quantum observables $C(p,q)$, $p,q\in {\mathbb Z}$. 
The equality
$S^{-1}C(p,q)S=C(-q,p)$ implies
\begin{eqnarray*}
C(p,q)S\zeta_j(z)=SC(-q,p)\zeta _j(z), \quad j=1,2,\ldots, r-1.
\end{eqnarray*}
At this point we need to make sure that we are able to shift indices
in the summation, and for that we have to let these indices range 
between $1$ and $2r$ (and not just between $1$ and $r-1$). To this end
we write $S\zeta_j(z)=\sum_{k=1}^{2r}\frac{1}{2}a_{k,j}\zeta_k(z)$,
$j=1,2,\ldots, 2r$. This is no longer an expansion in the basis of the 
Hilbert space, as each element of the basis appears twice.
Consequently, for a fixed $j$ we have 
\begin{eqnarray*}
&& \sum_{k=1}^{2r}\left(t^{-pq+2q(k+p)}a_{k+p,j}+t^{-pq-2q(k-p)}a_{k-p,j}
\right)\zeta _k(z)\\
&& = \sum_{k=1}^{2r}\left(t^{pq+2pj}a_{k,j+q}+t^{pq-2pj}a_{k,j-q}\right)
\zeta_k(z).
\end{eqnarray*}
Both sides of the equality are antisymmetric under
$k\rightarrow{}2r-k$.
For this reason we can equation the 
the coefficients of $\zeta_k(z)$ to obtain that for any $p,q,k,j$,
\begin{eqnarray*}
t^{2qk}a_{k+p,j}+t^{-2qk}a_{k-p,j}=t^{2pj}a_{k,j+q}+
t^{-2pj}a_{k,j-q}.
\end{eqnarray*}
Setting $p=0, q=k=1$ we obtain the recursive relation
\begin{eqnarray*}
a_{1,j+1}=(t^2+t^{-2})a_{1,j}-a_{1,j-1}.
\end{eqnarray*}
Since we are looking for a projective representation, we can set $a_{1,1}=1$,
which combined with $a_{1,0}=0$ yields $a_{1,j}=[j]$. 

Also, setting $q=0, p=1$, we obtain the recursive relation 
\begin{eqnarray*}
a_{k+1,j}+a_{k-1,j}=(t^{2j}+t^{-2j})a_{k,j},
\end{eqnarray*}
whence inductively we obtain $a_{k,j}=[kj]$, as desired.

Let us study the $T$-matrix now.  Similarly 
\begin{eqnarray*}
C(p,q)T\zeta _j(z)=TC(p,q+p)\zeta _j(z), \quad j=1,2,\ldots, r-1.
\end{eqnarray*} 
Again we extend the indices to the full range $1$ through $2r$ to be
able to shift indices in the summation, and write the above equality
in expanded form as 
\begin{eqnarray*}
&& \sum_{k=1}^{2r}\left(t^{-pq+2q(k+p)}b_{k+p,j}+
t^{-pq-2q(k-p)}b_{k-p,j}\right)\zeta _k(z)\\
&&  =\sum_{k=1}^{2r}\left(t^{-p(q+p)+2(q+p)j}b_{k,j-p}+t^{-p(q+p)-
2(q+p)j}b_{k,j+p}\right)\zeta_k(z)
\end{eqnarray*}
Hence for any $p,q,k,j$,
\begin{eqnarray*}
t^{2qk+2pq}b_{k+p,j}+t^{-2qk+2pq}b_{k-p,j}=t^{-p^2+2qj+2pj}b_{k,j-p}+t^{-p^2-
2qj-2pj}b_{k,j+p}.
\end{eqnarray*}
For $p=0$ we obtain 
\begin{eqnarray*}
(t^{2qk}+t^{-2qk})b_{k,j}=(t^{2qj}+t^{-2qj})b_{k,j}.
\end{eqnarray*}
This implies that $b_{k,j}=0$ if $k\neq j$, therefore $T$ is diagonal.
 Setting $p=1,k=j-1$ we obtain
\begin{eqnarray*}
b_{j,j}=t^{2j-1}b_{j-1,j-1}, \quad j=1,2,\ldots, r-1.
\end{eqnarray*}
Again, since we are looking for a projective representation, we
are allowed to choose $b_{1,1}=1$, in which case we obtain
inductively  $b_{j,j}=t^{j^2-1}$, $j=1,2,\ldots, r-1$, as desired. 
\end{proof}

We stress again that this theorem shows how the 
well known projective 
representation of the mapping class group of the torus 
on the Hilbert space can be introduced naturally using Weyl quantization. 

It is time now to describe the action of $S$ and $T$ on the vectors of
the basis. There is nothing to discuss about $T$ since  it is diagonal. 
For $S$ we  have

\begin{proposition}\label{Action}
The action of $S$ on the basis $\zeta _m(z)$, $m=1,2,\ldots, r-1$ is
given by
\begin{eqnarray*}
S\zeta_m(z)
=2i\sqrt{2}r^{\frac{3}{4}}e^{\frac{-\pi m^2}{2r}}\sum_{k=-\infty}^\infty
e^{-2\pi r(z-k)^2}\sinh 2\pi m (z-k)
\end{eqnarray*}
\end{proposition}

\begin{proof} We have $S\zeta _m(z)
=\frac{1}{X}\sum _{j=1}^{r-1}[jm]\zeta _j(z)$.
By Lemma \ref{formula} this is equal to 
\begin{eqnarray*}
\frac{2i\sqrt{2}}{\sqrt[4]{r}}\sin^2\frac{\pi}{r}
e^{-2\pi rz^2}\sum_{n=-\infty}^\infty
e^{2\pi nz-\frac{\pi n^2}{2r}}\sum _{j=1}^{r-1}[nj][jm].
\end{eqnarray*}

We compute 
$\sum_{j=1}^{r-1}[nj][jm]$, which  is  
\begin{eqnarray*}
&&(t^2-t^{-2})^{-2}\sum_{j=0}^{r-1}\left(e^{\frac{\pi inj}{r}}-
e^{-\frac{\pi inj}{r}}\right)\left(e^{\frac{\pi ijm}{r}}
-e^{-\frac{\pi ijm}{r}}\right)\\
& & \quad =
(t^2-t^{-2})^{-2}\sum_{-r+1\leq j\leq r-1}
\left(e^{\frac{\pi ij(n+m)}{r}}-e^{\frac{\pi ij(n-m)}{r}}\right).
\end{eqnarray*}
But
\begin{eqnarray*}
\sum_{-r+1\leq j\leq r-1} \left(e^{\frac{\pi ik}{r}}\right)^j
=
\left\{ 
\begin{array}{lc}
2r-1 & \mbox{if } 2r \mbox{ divides } k\\
-(-1)^{k}&  \mbox{otherwise}.
\end{array}
\right.
\end{eqnarray*}
Therefore the sum we are computing is
equal to  
\begin{eqnarray*}
2r\sum_{ n,2r|n+m}e^{-2\pi rz^2-\frac{\pi n^2}{2r}+2\pi nz}-
2r\sum_{ n,2r|n-m}e^{-2\pi rz^2-\frac{\pi n^2}{2r}+2\pi nz}
\end{eqnarray*}
multiplied by  $ir^{-\frac{1}{4}}/\sqrt{2}$.
Writing $n=2rk\pm m$ we obtain the formula from the statement.  
\end{proof}

\section{Knot and link  invariants as holomorphic sections}

In the Reshetikhin-Turaev theory, the quantum invariant of
a knot, viewed as a vector in the Hilbert space associated to the
torus, is expressed as
\begin{eqnarray*}
Z(K)=\frac{1}{X}\sum_{k=1}^{r-1}J(K,j)V^j(\alpha).
\end{eqnarray*}
where $X=\sqrt{\sum_{k=1}^{r-1}[k]^2}$, 
 $J(K,j)$ is the $j$th colored Jones polynomial of $K$,
and $V^j(\alpha)$ is the orthonormal basis consisting of 
colorings of the core of the solid torus by irreducible representations.
As explained in Section 2, the result from \cite{GU} allows us to 
identify the elements of this orthonormal basis with holomorphic
sections of the Chern-Simons line bundle. We therefore have

\begin{proposition}\label{invariant}
The quantum invariant in level $r$ of a knot $K$ is the holomorphic section
of the Chern-Simons line bundle over ${\mathcal M}$ defined by the formula
\begin{eqnarray*}
Z(K)=2\sqrt{2}ir^{-\frac{1}{4}}\sin^2\frac{\pi}{r}
\sum_{n=-\infty}^\infty
e^{-2\pi r\left(z-\frac{ n}{2r}\right)^2}\sum _{j=1}^{r-1}[nj]J(K,j),
\end{eqnarray*}
where $J(K,j)$ is the $j$th colored Jones polynomial of $K$.
\end{proposition}

\begin{proof} It follows from Lemma \ref{formula} since 
\begin{eqnarray*}
Z(K)=\frac{1}{X}\sum_{j=1}^{r-1}J(K,j)\zeta _j(z).\qedhere
\end{eqnarray*}
\end{proof}

\begin{example}\label{unknot}
The quantum invariant of the trivial knot is 
\begin{eqnarray*}
Z(0)=2i\sqrt{2}r^{\frac{3}{4}}e^{-\frac{\pi}{2r}}\sum_{n=-\infty}^\infty
e^{-2\pi r(z-n)^2}\sinh 2\pi (z-n).
\end{eqnarray*}
\end{example}

Because $J(0,j)=[j]$, 
the formula is a particular case of  Proposition \ref{Action}, since
$Z(0)=S\zeta _1(z)$.



\begin{example}
The quantum invariant 
of the $(p,q)$-torus knot is 
\begin{eqnarray*}
Z(K_{p,q})=-\frac{1}{\sqrt{2}}r^{-\frac{1}{4}}\sin\frac{\pi}{r}
\sum_{n=-\infty}^\infty C_n
e^{-2\pi r\left( z-\frac{ n}{2r}\right)^2},
\end{eqnarray*}
where 
\begin{eqnarray*}
C_n=\frac{1}{\sin \frac{n\pi}{r}}\sum_{k=1}^{r-1}
t^{-pqk^2}\left(\left[2n\left\lfloor \frac{r-1-k}{2}\right\rfloor +kn+n\right]
-[kn-n]\right)\\\times
\left([kp+kq+1]-[kp-kq+1]\right).
\end{eqnarray*}
In this formula  square brackets represent quantized integers while
$\lfloor \cdot \rfloor $ represents the greatest integer function.
\end{example}

This is a consequence of the formula for the $j$th colored 
Jones polynomial of a torus knot
deduced in \cite{gelca1}:
\begin{eqnarray*}
J(K_{p,q},j)=\sum_{\begin{array}{c}
0\leq k\leq j\\ k\equiv j(\mbox{mod }2)
\end{array}}\frac{t^{-pqk^2}}{t^2-t^{-2}}([kp+kq+1]-[kp-kq+1]).
\end{eqnarray*}
Explicitly
\begin{alignat*}{1}
Z(K_{p,q})&=2\sqrt{2}ir^{-\frac{1}{4}}\sin^2\frac{\pi}{r}
\sum_{n=-\infty}^\infty
e^{-2\pi r\left(z-\frac{ n}{2r}\right)^2}\sum
_{j=1}^{r-1}[nj]J(T_{p,q},j)\\
&=-\frac{i}{\sqrt{2}}r^{-\frac{1}{4}}\sin\frac{\pi}{r}
\sum_{n=-\infty}^\infty 
e^{-2\pi r\left( z-\frac{ n}{2r}\right)^2}\\
& \quad\times \sum_{j=1}^{r-1}[nj]\sum_{\begin{array}{c}
0\leq k\leq j\\ k\equiv j(\mbox{mod }2)
\end{array}}t^{-pqk^2}([kp+kq+1]-[kp-kq+1]).
\end{alignat*}
Changing the order of summation in the double sum we find that it is
equal to
\begin{eqnarray*}
&&\sum_{k=1}^{r-1}t^{-pqk^2}([kp+kq+1]-[kp-kq+1])\sum_{\begin{array}{c}
k\leq j\leq r-1\\ j\equiv k(\mbox{mod }2)\end{array}}[nj]\\
&& \quad =\sum_{k=1}^{r-1}t^{-pqk^2}([kp+kq+1]-[kp-kq+1])\sum_{
0\leq m\leq \lfloor \frac{r-1-k}{2}\rfloor}[n(2m+k)].
\end{eqnarray*}
Write the quantized integer in explicit form, then sum the
exponentials as geometric series to obtain that the inside sum
is equal to 
\begin{eqnarray*}
&&e^{-\frac{n\pi i}{r}}\frac{\exp \left(\frac{n\pi i}{r}\left(2\left(
\left\lfloor \frac{r-1-k}{2}\right\rfloor +1\right) +k\right)\right)-
\exp \left(\frac{nk\pi i}{r}\right)}{\exp\left(\frac{n\pi i}{r}\right)-\exp
  \left(-\frac{n\pi i}{r}\right)}\\
&& \quad -e^{\frac{n\pi i}{r}}\frac{\exp \left(-\frac{n\pi i}{r}\left(2\left(
\left\lfloor \frac{r-1-k}{2}\right\rfloor +1\right) +k\right)\right)-
\exp \left(-\frac{nk\pi i}{r}\right)}{\exp\left(\frac{n\pi i}{r}\right)-\exp
  \left(-\frac{n\pi i}{r}\right)}
\end{eqnarray*}
everything multiplied by a factor of $(e^{\frac{\pi
    i}{r}}-e^{-\frac{-\pi i}{r}})^{-1}$. Using the definition of 
quantized integers, we find that this is equal to $C_n$. 

Because  of the presence of the greatest integer function,
the formula  cannot be further simplified using a  Gauss sum. 

For a link $L$ with $k$ components, the Hilbert space is obtained
by taking the tensor product of $k$ copies of the Hilbert space of  the torus.
The formula for the quantum invariant is then
\begin{eqnarray*}
\frac{1}{X}
\sum_{j_1,j_2,\cdots,j_k=1}^{r-1}J(L,j_1,j_2,\cdots, j_k)V^{j_1}(\alpha)
\otimes V^{j_2}(\alpha)\otimes \cdots \otimes V^{j_k}(\alpha).
\end{eqnarray*}
This can again be translated to the analytical setting 
replacing $V^{j}(\alpha)$'s by $\zeta_{j}(z)$'s. Here is one
example. 

\begin{example}
The quantum invariant 
of the Hopf link is
\begin{eqnarray*}
Z(L)=-r\sqrt{2}\sin \frac{\pi}{r}\sum_{n=-\infty}^\infty
\sum_{k=-\infty}^\infty [nk]e^{-2r\pi \left[ \left(z-\frac{n}{2r}\right)^2+
\left(w-\frac{k}{2r}\right)^2\right]}.
\end{eqnarray*}
\end{example}

Indeed, 
\begin{eqnarray*}
Z(L)=\frac{1}{X}\sum_{j,m=1}^{r-1}[jm]\zeta_j(z)\zeta _m(w).
\end{eqnarray*}
Using Lemma \ref{formula} we transform this into
\begin{eqnarray*}
-\frac{i}{2\sqrt{2}}\sum_{n=-\infty}^\infty
\sum_{k=-\infty}^\infty e^{-2r\pi \left[ \left(z-\frac{n}{2r}\right)^2+
\left(w-\frac{k}{2r}\right)^2\right]}\sum_{m,j=1}^{r-1}(e^{\frac{\pi inj}{r}}-
e^{-\frac{\pi inj}{r}})\times\\
\times (e^{\frac{\pi ijm}{r}}-
e^{-\frac{\pi ijm}{r}})(e^{\frac{\pi imk}{r}}-
e^{-\frac{\pi imk}{r}}).
\end{eqnarray*}
A computation with roots of unity similar to the one from the proof of 
Proposition \ref{Action} shows that the innermost double sum is equal
to 
\begin{eqnarray*}
-2r(e^{\frac{\pi ink}{r}}-e^{-\frac{\pi ink}{r}})
\end{eqnarray*}
 and the formula follows.

\section{Some properties of the quantum observables}

In the Feynman path integral formulation, the operator $C(p,q)$ 
representing the quantization
of the function $2\cos 2\pi (px+qy) $ is  the integral operator with kernel
\begin{eqnarray*}
{\mathcal K}_{p,q}(A_1,A_2)=\int_{{\mathcal M}_{A_1,A_2}} e^{iN{\mathcal L}
(A)}(\mbox{tr}_{V^{n+1}}-\mbox{tr}_{V^{n-1}}) (\mbox{hol}_C(A)) {\mathcal D}A.
\end{eqnarray*}
Here 
$A_1$, $A_2$ are conjugacy classes of  
connections on the torus ${\mathbb T}^2$ modulo gauge transformations, 
$A$ is a conjugacy class
of connections on  ${\mathbb T}^2\times [0,1]$ modulo gauge transformations
such that $A|_{{\mathbb T}^2\times\{0\}}=A_1$
and $A|_{{\mathbb T}^2\times \{1\}}=A_2$, $n$ is the greatest common
divisor of $p$ and $q$,  and $\mbox{tr}_{V^n}
(\mbox{hol}_C(A))$, known as the Wilson line, is 
the trace of the $n$-dimensional irreducible representation
of $SU(2)$ evaluated on the  holonomy
of $A$ around the curve
 $C$ of slope $p/q$. The ``integral'' is taken over all 
conjugacy classes of  connections $A$. 

We now exhibit a mathematically 
well defined formula for this kernel. 
In complex coordinates, the kernel  
 is given by 
\begin{eqnarray*}
{\mathcal K}_{p,q}(z,w)=\sum_{j=1}^{r-1} \left(C(p,q)\zeta_j\right)(z)
\overline{\zeta _j(w)}.
\end{eqnarray*}
 A straightforward computation 
shows that 

\begin{proposition}
The kernel of the operator $C(p,q)$ is given by 
\begin{alignat*}{1}
{\mathcal K}_{p,q}(z,w)=&2r^{\frac{3}{2}}\sum_{\begin{array}{l}
m,n=-\infty\\
2r|q\pm (n-m)
\end{array}
}^\infty e^{-2r\pi \left[\left(z-\frac{n}{2r}\right)^2+\left(\overline{w}
-\frac{
m}{2r}\right)^2\right]\mp \frac{npi \pi }{r}}\\
&-2r^{\frac{3}{2}}\sum_{\begin{array}{l}
m,n=-\infty\\
2r|q\pm (n+m)
\end{array}
}^\infty e^{-2r\pi \left[\left(z-\frac{n}{2r}\right)^2+\left(\overline{w}
-\frac{
m}{2r}\right)^2\right]\mp \frac{npi \pi }{r}}.
\end{alignat*}
\end{proposition}

The operator $C(p,q)$ acts on theta functions in the Hilbert space
${\mathcal H}_r$ by
\begin{eqnarray*}
(C(p,q)f)(z)=\int _{\mathcal M}{\mathcal K}_{p,q}(z,w)f(w)dw.
\end{eqnarray*}
Clearly $K_{p,q}(z,w)$ is holomorphic in $z$ and antiholomorphic in $w$. 


\begin{proposition} The characteristic polynomial of the
operator $C(p,q)$ is 
\begin{eqnarray*}
\prod_{k=1}^{r-1}\left(\lambda -2\cos
\frac{\mbox{gcd}(p,q,2r)k\pi}{r}\right).
\end{eqnarray*}
\end{proposition}

\begin{proof}
Note that if $p=np'$, $q=nq'$, with $p',q'$ coprime, then
$C(p,q)=T_{n}(C(p',q'))$, where $T_n(x)$ is the $n$th Chebyshev
polynomial
(subject to the normalization $T_0(x)=2$, $T_1(x)=x$, $T_{n+1}(x)=
xT_n(x)-T_{n-1}(x)$, $n\geq 1$).
So the case where $p$ and $q$ have a common divisor follows from
the case where they are coprime via the spectral mapping theorem. 
Let us assume that  $p$ and $q$ are coprime. 
As a consequence of Theorem \ref{unique}, there exists an invertible  matrix
$A$ such that $C(p,q)=A^{-1}C(1,0)A$. In fact  $A=\rho(g)$ where $g$ is
the element of the mapping class group that maps the $(1,0)$ curve on
the torus to the $(p,q)$ curve. Because the characteristic polynomial
is invariant under conjugation, it suffices to prove the property for 
$C(1,0)$. It is easy to check that the characteristic polynomials
in dimensions $r+1$, $r$, and $r-1$ are related by $p_{r+1}(\lambda)=\lambda
p_r(\lambda)-p_{r-1}(\lambda)$. Also $p_1(\lambda )=\lambda$ and $p_2(\lambda)=
\lambda^2-1$. It follows that $p_{r}(\lambda)=S_{r}(\lambda)$,
 where $S_n(\lambda)$ denotes the
Chebyshev polynomial of second type, $S_0(\lambda)=1$, $S_1(\lambda )=\lambda$,
$S_{n+1}(\lambda)=\lambda S_n(\lambda)-S_{n-1}(\lambda)$, $n\geq 1$.
Factoring $S_r(\lambda)$ we obtain the desired formula. 
\end{proof}

As a corollary, we recover in analytical setting the well known
topological  fact that a simple closed
curve on the torus colored
by the $r$-dimensional irreducible representation of the quantum
group of $SL(2,{\mathbb C})$ is equal to zero, i.e. the
operator with symbol equal to the trace of the holonomy along a curve in the
$r$-dimensional irreducible representation of $SU(2)$ is the zero operator. 

\section{An application to quantum computing}

In \cite{fnsww} the authors analyzed a possible quantum system suitable
for quantum computation which is based on the fractional quantum 
Hall system. The model we have in mind happens at the plateau
corresponding to the fraction $12/5$, where a non-abelian statistics
has been predicted. 
The subspace of  ground states of the Hilbert space of the
quantum system can be identified with the vector space of an 
$SU(2)_r\times \overline{SU(2)}_r$ Chern-Simons quantum field theory
for $r=5$.
This in turn can be obtained through Drinfeld's double construction, or
can be simply identified with the linear space of operators 
(quantum observables) of the $SU(2)_r$ Chern-Simons theory. 
The authors considered the case of the torus and were particularly 
interested in finding a basis of this vector space in terms of 
curves on the torus colored by representations of the quantum group
of $SU(2)$. They succeeded for the case  $r=3$, but the real goal was 
 $r=5$. In this section we will
answer their question for an arbitrary $r$. 
The problem was bought to our attention by Zh. Wang.

It is known that the vector space in discussion has dimension
$(r-1)^2$ and is generated by the operators $C(p,q)$, $p,q\in {\mathbb Z}$
introduced  before. Recall that if $n$ denotes the greatest common divisor
of $p$ and $q$, and $p'=p/n$, $q'=q/n$, then $C(p,q)$ can be identified
with the curve of slope $p'/q'$ on the torus colored by the difference
of the $n+1$st and $n-1$st dimensional irreducible representations
of $SU(2)$. Our goal is to find a basis consisting of operators of
the form $C(p,q)$. The key idea is to work in the basis $\zeta_{j}(z)$,
$j=1,2,\cdots, r-1$, and use the formula
\begin{eqnarray*}
C(p,q)\zeta_m(z)=t^{-pq}(t^{2qm}\zeta_{m-p}(z)+t^{-2qm}\zeta_{m+p}(z)).
\end{eqnarray*}
We see now that in the basis 
$\zeta_j(z)$, $j=1,2,\cdots, r-1$ the matrices of the
operators $C(p,q)$ are particularly simple. 
To summarize our approach, through a linear isomorphism 
we identify the Hilbert space with a space of operators, then choose a basis
in which the matrices of these operators are simple enough. We obtain 

\begin{theorem}
A basis of the linear space of quantum observables of the $SU(2)_r$
Chern-Simons theory on the torus is given by the operators 
\begin{eqnarray*}
&& C(0,q), \quad 0\leq q\leq r-2,\\
&& C(p,q), \quad 1\leq p\leq r-2, -r+p+2\leq q\leq r-p-1.
\end{eqnarray*}
\end{theorem}

\medskip

\begin{proof} 
We will 
show that the $C(p,q)$, with $p,q$ ranging as described 
in the statement, span the space 
of $(r-1)\times (r-1)$ matrices. Start with the diagonal. 

\begin{lemma}
The diagonal matrices are  spanned by  $C(0,q)$, $0\leq q\leq r-2$.
\end{lemma}

\begin{proof}
 The matrix of $C(0,q)$ is the diagonal matrix with
entries $$\left(\cos \frac{q\pi}{r}, \cos\frac{2q\pi}{r},
\ldots,\cos\frac{(r-1)q\pi}{r}\right).$$ 
Denote $\alpha =\frac{\pi}{r}$. We
have to show that the determinant  
\begin{eqnarray*}
\left|
\begin{array}{llll}
1& 1&\cdots &1\\
\cos \alpha &\cos 2\alpha & \cdots & \cos (r-1)\alpha\\
\cos 2\alpha &\cos 4\alpha &\cdots &\cos 2(r-1)\alpha\\
\cdots&\cdots &\ddots &\cdots\\
\cos(r-2)\alpha&\cos 2(r-2)\alpha &\cdots &\cos (r-1)(r-2)\alpha
\end{array}
\right|
\end{eqnarray*}
is nonzero. 
To compute the determinant, let $x_1=2\cos \alpha$, $x_2=2\cos 2\alpha$,
$\ldots $, $x_{r-1}=2\cos (r-1)\alpha$. Denote by $T_n(x)$ the $n$th
Chebyshev polynomial ($T_0(x)=2$, $T_1(x)=x$, $T_{n+1}(x)=xT_n(x)-T_{n-1}(x))$.
Then the determinant is
\begin{eqnarray*}\frac{1}{2^{r-1}}
\left|
\begin{array}{llll}
T_0(x_1)& T_0(x_2)&\cdots &T_0(x_{r-1})\\
T_1(x_1)&T_1(x_2)&\cdots &T_1(x_{r-1})\\
\cdots &\cdots &\ddots &\cdots\\
T_{r-2}(x_1)&T_{r-2}(x_2)&\cdots &T_{r-2}(x_{r-1})
\end{array}
\right|
\end{eqnarray*}
Row operations transform this into the Vandermonde determinant. We
conclude that  the
value of the original determinant is
\begin{eqnarray*}
\frac{1}{2^{r-1}}
\prod_{1\leq k<j\leq r-1}\left(\cos \frac{j\pi}{r}-\cos \frac{k\pi}{r}\right)
\neq 0.\qedhere
\end{eqnarray*}
\end{proof}

Let us return to the proof of the theorem.  
 For some nonzero $k\leq r-2$, let us look at those 
 $C(p,q)$  $ 0\leq p\leq k$. The nonzero entries  of the  matrix 
of such an element lie at distance at most $k$ from the main
diagonal, i.e. they are among the  $a_{ij}$'s with $i-k\leq j\leq i+k$. 

We prove by induction on $p$ that $C(k,q)$ with $q$ subject to the conditions
from the statement and  $0\leq k\leq p$
span $M_p$, the set of  all matrices whose only nonzero elements are
of the form $a_{ij}$, with $i-p\leq j\leq i+p$.
The base case $p=0$ was proved in the lemma. 

Assume that the property is
true for $p-1$, and let us prove it for $p$.
Consider the matrix of $t^{pq}C(p,q)$. Using the inductive hypothesis
we can  add to it an element
of $M_{p-1}$, so that the resulting matrix  
$A_{p,q}$ is of the form
\begin{eqnarray*}
\left(
\begin{array}{lllllll}
0&0&\cdots &t^{-2q}&0&\cdots &0\\
0&0&\cdots &0&t^{-4q}&\cdots &0\\
\cdots &\cdots &\ddots &\cdots &\cdots &\cdots &\cdots\\
0&0&\cdots &0&0&\cdots & t^{-2(r-1-p)q}\\
\cdots &\cdots &\ddots &\cdots &\cdots &\ddots&\cdots\\
t^{2(p+1)q}&0&\ldots &0&0&\ldots &0\\
0&t^{2(p+2)q}&\ldots &0&0&\ldots &0\\
\cdots &\cdots &\ddots &\cdots &\cdots &\ddots &\cdots\\
\end{array}
\right).
\end{eqnarray*}

The nonzero entries
of $A_{p,q}$ are those of indices $(1,p+1)$, $(2,p+2)$, $\ldots$,
$(r-p-1,r-1)$, $(p+1,1)$, $(p+2,2)$,$\ldots $, $(r-1,r-p-1)$ (those
at distance $p$ from the main diagonal of the matrix).
The space $M_p\ominus M_{p-1}$ has dimension $2r-2k-2$, and
a basis $E_{i,j}$, $(i,j)\in \{(1,p+1),(2,p+2), \ldots\}$, where
$E_{i,j}$ denotes the matrix whose only nonzero entry is equal to $1$ and
is that of index $(i,j)$. In this basis the coordinates of $A_{p,q}$
are 
\begin{eqnarray*}
\left(t^{(-2)q}, t^{(-4)q},\ldots, t^{-2(r-p-1)q},t^{2(p+1)q},\ldots, 
t^{2(r-1)q}\right).
\end{eqnarray*}
To show that $A_{p,q}$, $-r+p+2\leq q\leq r-p-1$
form a basis of $M_{p}\ominus M_{p-1}$ we arrange the entries of these vectors
 in a determinant, and show that this determinant is not equal to $0$.
With the convention $x_1=t^{-2},x_2=t^{-4}, x_3=t^{-6}$, ..., 
the determinant is 
\begin{eqnarray*}
\left|
\begin{array}{rrrr}
x_1^{-r+p+2} & x_2^{-r+p+2}&\ldots &x_{2r-2p-2}^{-r+p+2}\\
x_1^{-r+p+3}& x_2^{-r+p+3}&\ldots &x_{2r-2p-2}^{-r+p+3}\\
\cdots &\cdots &\ddots &\cdots \\
x_1^{r-p-1}&x_2^{r-p-1}&\ldots &x_{2r-2p-2}^{r-p-1}
\end{array}
\right|.
\end{eqnarray*}
Multiplying this determinant by
\begin{eqnarray*}
(x_1x_2\cdots x_{2r-2p})^{r-p-2}
\end{eqnarray*}
produces  a Vandermonde determinant, which is nonzero since the $x_i$'s
are distinct. This completes the inductive argument, and consequently the proof
of the theorem.
\end{proof}

We point out that contrary to a naive intuition, the indices of the basis
elements do not range in an $(r-1)\times (r-1)$ rectangle, but in a 
triangular region, a surprising 
fact already observed in \cite{fnsww} for $r=3$. 
For $r=5$, we would like to describe a basis more 
in the spirit of the above mentioned
paper. For that let us denote
 by $V^k(m,n)$ the curve of slope $m/n$ on the torus 
colored by the $k$-dimensional irreducible representation of the
quantum group of $SL(2,{\mathbb C})$. 

\begin{cor}
The  linear space of quantum observables for $r=5$ has a 
basis  formed by the identity operator together with
the operators 
$V^2(0,1)$, $V^3(0,1)$, $V^4(0,1)$, $V^2(1,-2)$, $V^2(1,-1)$,
$V^2(1,0)$, $V^2(1,1)$, $V^2(1,2)$,  $V^2(1,3)$, 
$V^2(2,-1)$, $V^3(1,0)$, $V^2(2,1)$, $V^3(1,1)$, $V^4(1,0)$, $V^2(3,1)$. 
\end{cor}

\begin{proof}
This follows from the theorem using the identity
\begin{eqnarray*} 
C(p,q)=V^{n+1}(p',q')-V^{n-1}(p',q'),
\end{eqnarray*}
where $n$ is the greatest common divisor of $p$ and $q$, and $p'=p/n$,
$q'=q/n$. 
\end{proof}

With the usual conventions for curves (for example that $(4,2)$ means
the double of the curve $(2,1)$), we can rephrase this as

\begin{cor}
The  linear space of quantum observables for $r=5$ has a 
basis  formed by the identity operator together with
the operators 
$(0,1)$, $(0,2)$, $(0,3)$, $(1,-2)$, $(1,-1)$,
$(1,0)$, $(1,1)$, $(1,2)$,  $(1,3)$, 
$(2,-1)$, $(2,0)$, $(2,1)$, $(2,2)$, $(3,0)$, $(3,1)$. 
\end{cor}


\begin{thebibliography}{000}


\bibitem{jea} J.E. Andersen, {\em Deformation quantization and
 geometric quantization
of Abelian moduli spaces}, Comm. Math. Phys. {\bf 255} (2005), no. 3,
727--745.

\bibitem{jeau1}
J.E. Andersen, K. Ueno, {\em Geometric construction of modular functors from
conformal field theory}, preprint.

\bibitem{jeau2}
J.E. Andersen, K. Ueno, {\em Abelian conformal field theories and
determinant bundles}, preprint.

\bibitem{folland} G. Folland, {\em Harmonic Analysis in Phase
  Space}, Princeton University Press, 1989.

\bibitem{freed}
D. Freed, {\em Classical Chern-Simons}, part 1 Adv. Math. {\bf 113}
(1995) no. 2, 237--303, part II Houston J. Math. {\bf 28} (2002)
no. 2, 293--310.

\bibitem{fnsww} M. Freedman, Ch.  Nayak, K. Shtengel, 
K. Walker, Zh.  Wang, {\em A class of $P,T$-invariant topological
phases of interacting electrons}, Annals of Physics, {\bf 310}(2004), 
428--492.


\bibitem{GU} R. Gelca, A. Uribe, {\em The  Weyl 
quantization and the quantum group quantization of the moduli 
space of flat $SU(2)$-connections on the torus are the same},
Commun. Math. Phys., {\bf 233}(2003), 493--512. 

\bibitem{gelca1} R. Gelca, {\em The quantum invariant of the
complement of a regular neighborhood of a link}, Topology and its Appl.,
{\bf 81}(1997, 147--157.


\bibitem{Jones} V.F.R. Jones, {\em Polynomial invariants of knots via von
Neumann algebras}, 
Bull. Amer. Math. Soc., {\bf 12}(1985), 103--111.

\bibitem{RT} N.Yu. Reshetikhin,  V.G. Turaev,
 {\em Invariants of 3-manifolds
via link polynomials and quantum groups}, { Inventiones Math.},
 {\bf 103}(1991),
547--597.


\bibitem{T2} V.G. Turaev, {\em Quantum Invariants of Knots and 
3-Manifolds}, de Gruyter, 1994.

\bibitem{Wi} E. Witten, {\em Quantum field theory and 
the Jones polynomial},
Comm. Math. Phys., {\bf 121}(1989), 351--399.

\end{thebibliography}
\end{document}